\documentclass[11pt]{amsart}

\usepackage{amsmath,amssymb,amsthm}
\usepackage[a4paper,margin=25mm]{geometry}
\usepackage{booktabs}
\usepackage{xcolor}
\usepackage{tikz-cd}
\usepackage{hyperref}

\hypersetup{
  colorlinks=true,
  linkcolor=blue!45!black,
  citecolor=blue!45!black,
  urlcolor=blue!45!black
}

\theoremstyle{plain}
\newtheorem{theorem}{Theorem}[section]
\newtheorem{proposition}[theorem]{Proposition}
\newtheorem{lemma}[theorem]{Lemma}
\newtheorem{corollary}[theorem]{Corollary}

\theoremstyle{definition}

\newtheorem{example}[theorem]{Example}

\theoremstyle{remark}
\newtheorem{remark}[theorem]{Remark}

\theoremstyle{definition}
\newtheorem*{ack}{Acknowledgement}

\newcommand{\ZZ}{\mathbb Z}
\newcommand{\cI}{\mathcal I}
\newcommand{\cO}{\mathcal O}
\newcommand{\cB}{\mathcal B}
\newcommand{\cC}{\mathcal C}
\newcommand{\setn}[1]{\{1,\ldots,#1\}}
\newcommand{\class}[1]{\left[#1\right]_{\operatorname{H}}}
\DeclareMathOperator{\Ker}{Ker}

\begin{document}

\title{Explicit Bases for the Second Integral Homology of Cactus Groups}

\author{JINGYA ZHANG}

\address{Department of Mathematics, Graduate School of Science, Hokkaido University, Sapporo, 060-0810 Japan}

\email{jingya.zhang.y6@elms.hokudai.ac.jp}


\begin{abstract}
In this paper, we determine the second integral homology $H_2(J_n;\ZZ)$ of cactus groups $J_n$.
We first use Hopf's formula and Pontryagin products to obtain a smaller generating set of $H_2(J_n;\ZZ)$.
We then construct a detecting homomorphism to prove that these generators are linearly independent and hence form an explicit basis.
\end{abstract}

\maketitle

\section{Introduction}

Cactus groups $J_n\,(n\geq 2)$ originated as \emph{quasibraid groups} in Devadoss~\cite{Devadoss} and \emph{mock reflection groups} in Davis, Januszkiewicz, and Scott~\cite{DavisJanuszkiewiczScott}.
Henriques and Kamnitzer introduced the name \emph{cactus group} in the context of coboundary categories and crystals~\cite{HenriquesKamnitzer}.
Let $\overline{\mathcal M}_{0,n+1}(\mathbb R)$ denote the real locus of the moduli space of stable genus-zero curves with $n+1$ marked points.
The symmetric group $S_n$ acts on this space by permuting the first $n$ marked points.
The cactus group $J_n$ is isomorphic to the orbifold fundamental group of the quotient orbifold $[\overline{\mathcal M}_{0,n+1}(\mathbb R)/S_n]$.

Concretely, for $1\leq p<q\leq n$, let $s_{p,q}$ be the generator associated with the interval $[p,q]=\{p,p+1,\ldots,q\}$.
The \emph{cactus group} $J_n\,(n\geq 2)$ is defined by these generators and the following relations:
\begin{itemize}
	\item $s_{p,q}^2=1,\qquad 1 \leq p<q \leq n$,
	\item
$s_{p,q}s_{r,t}=s_{r,t}s_{p,q},\qquad [p,q]\cap[r,t]=\varnothing$,
	\item $s_{p,q}s_{r,t}s_{p,q}=s_{p+q-t,p+q-r},\qquad [r,t]\subset[p,q]$.
\end{itemize}

The kernel of the natural homomorphism $J_n\to S_n$, which sends $s_{p,q}$ to the reversal of $[p,q]$, is called the \emph{pure cactus group}.
It is isomorphic to the fundamental group of the aspherical moduli space $\overline{\mathcal M}_{0,n+1}(\mathbb R)$~\cite{DavisJanuszkiewiczScott}, whose rational cohomology ring was determined in~\cite{EtingofHenriquesKamnitzerRains}.
Recent results on real moduli spaces and pure cactus groups can be found in~\cite{KhoroshkinWillwacher}.
Since $\overline{\mathcal M}_{0,n+1}(\mathbb R)$ is aspherical and $J_n$ is the orbifold fundamental group of $[\overline{\mathcal M}_{0,n+1}(\mathbb R)/S_n]$, the group homology of $J_n$ is isomorphic to the orbifold homology of  $[\overline{\mathcal M}_{0,n+1}(\mathbb R)/S_n]$.

Cactus groups also act on crystals in coboundary categories~\cite{HenriquesKamnitzer} and arise in the study of the octahedron recurrence, Coxeter groups, cells, and related combinatorics~\cite{HenriquesKamnitzerOctahedron,Bonnafe,Losev,ChmutovGlickPylyavskyy,BellingeriCheminLebed}.
These connections place cactus groups at an intersection of topology, algebra, and combinatorics.

Throughout this paper, the commutator of two elements $g,h$ of a group $G$ is defined by
\[
[g,h]:=ghg^{-1}h^{-1}.
\]

For $2\leq k\leq n$, define the standard generators by $\sigma_k:=s_{1,k}$.
Let $F_n=F(\sigma_2,\ldots,\sigma_n)$ be the free group.
Chemin and Nanda~\cite[Theorem~A]{CheminNanda} gave a minimal presentation of $J_n$ on the generators $\sigma_2,\ldots,\sigma_n$. 
Zhang~\cite[Corollary~2.8]{Zhang} simplified the relators of this minimal presentation, obtaining the following three families:
\[
\begin{aligned}
&A_n=\{\sigma_i^2 \mid 2\leq i \leq n\},\\     
&B_n=\{b(i,j,k)=[\sigma_k \sigma_i \sigma_k, \sigma_j] \mid 2\leq i\leq j,\ i+j\leq k\leq n\},\\
&C_n=\{c(i,j,k)=[\sigma_{k-i}\sigma_k\sigma_{i+j},\sigma_j] \mid  i\geq 1,\ j\geq 2,\ i+j<k\leq n,\  i+j\leq k-i\}.    
\end{aligned}
\]
Let $N_n$ be their normal closure in $F_n$. 
Then $J_n=F_n/N_n$.

By Hopf's formula~\cite[Section II.5]{Brown}, we have
\begin{equation}
\label{eq:hopf}
 H_2(J_n;\ZZ)
 \cong
\frac{N_n\cap[F_n,F_n]}{[N_n,F_n]}.
\end{equation}
For $r\in N_n\cap[F_n,F_n]$, we denote its Hopf class by
\[
\class{r}:=
r \pmod{[N_n,F_n]}.
\]

We begin with the following consequence of~\cite[Theorem~1.1 and Section~3]{Zhang}.

\begin{theorem}
\label{thm:zhang}
$H_2(J_2;\ZZ)$ and $H_2(J_3;\ZZ)$ are trivial.
For $n\geq 4$, $H_2(J_n;\ZZ)$ is a finite elementary abelian $2$-group generated by the Hopf classes of the relators in $B_n\cup C_n$.
\end{theorem}

In view of Theorem~\ref{thm:zhang}, we henceforth regard $H_2(J_n;\ZZ)$ as an $\mathbb F_2$-vector space and use additive notation for Hopf classes.
In Section~\ref{subsec:pontryagin}, we give a proof of Theorem~\ref{thm:zhang} using Hopf's formula and Pontryagin products.

The second integral homology $H_2(J_n;\ZZ)$ is the Schur multiplier, which is related to central extensions and projective representations.
With the trivial $J_n$-action on $\mathbb C^\times$, the universal coefficient theorem and the divisibility of $\mathbb C^\times$ give
\[
H^2(J_n;\mathbb C^\times)\cong
\operatorname{Hom}_{\ZZ}(H_2(J_n;\ZZ),\mathbb C^\times).
\]
By Theorem~\ref{thm:zhang}, $H_2(J_n;\mathbb Q)=0$ for all $n\geq 2$.
Hence an explicit $\mathbb F_2$-basis describes the $2$-torsion classes invisible to rational homology.

Although Theorem~\ref{thm:zhang} provides a generating set for $H_2(J_n;\ZZ)$ for each $n\geq 4$, it does not determine either its exact dimension over $\mathbb F_2$ or an explicit basis.
The purpose of this paper is to determine both.

In $H_2(J_n;\ZZ)$, define
\[
\begin{aligned}
\cB_n&:=\{\beta_{i,j}:=\class{b(i,j,i+j)}=\class{[\sigma_{i+j}\sigma_i\sigma_{i+j},\sigma_j]} \mid 2\leq i\leq j,\ i+j\leq n\},\\
\cC_n&:=\{\gamma_{j,k}:=\class{c(1,j,k)}= \class{[\sigma_{k-1}\sigma_k\sigma_{j+1},\sigma_j]} \mid j\geq 2,\ j+2\leq k\leq n\}.  
\end{aligned}
\]

Our main result is as follows.

\begin{theorem}
\label{thm:main}
For $n\geq 4$, $H_2(J_n;\ZZ)$ is an $\mathbb F_2$-vector space with basis
$\cB_n\cup\cC_n$.
Consequently,
\[
H_2(J_n;\ZZ)\cong
\begin{cases}
0, & n=2,3,\\[2mm]
(\ZZ/2)^{|\cB_n|+|\cC_n|}, & n\geq 4,
\end{cases}
\]
where
\[
|\cB_n|+|\cC_n|=
\begin{cases}
\dfrac{(n-2)(3n-8)}{4},
& \text{$n$ is even},\\[2mm]
\dfrac{(n-3)(3n-5)}{4},
& \text{$n$ is odd}.
\end{cases}
\]
\end{theorem}

As an immediate consequence of Theorem~\ref{thm:main}, we obtain
\[
H^2(J_n;\mathbb C^\times)\cong (\ZZ/2)^{|\cB_n|+|\cC_n|}
\qquad (n\geq 4).
\]

In particular, the dimensions over $\mathbb F_2$ given by Theorem~\ref{thm:main} for $4\leq n\leq 10$ are as follows.

\begin{center}
\begin{tabular}{c|rrrrrrr}
\toprule
$n$&4&5&6&7&8&9&10\\
\midrule
$|\cB_n|$&1&2&4&6&9&12&16\\
$|\cC_n|$&1&3&6&10&15&21&28\\
$|\cB_n|+|\cC_n|$&2&5&10&16&24&33&44\\
\bottomrule
\end{tabular}
\end{center}
Moreover, $H_2(J_4;\ZZ)$, $H_2(J_5;\ZZ)$, and $H_2(J_6;\ZZ)$ have the following $\mathbb F_2$-bases, respectively:
\[
\begin{aligned}
  &\{\beta_{2,2},
  \gamma_{2,4}\},\\
  &\{\beta_{2,2},
  \beta_{2,3},
  \gamma_{2,4},
  \gamma_{2,5},
  \gamma_{3,5}\},\\
  &\{\beta_{2,2},
\beta_{2,3},
\beta_{2,4},
\beta_{3,3},
\gamma_{2,4},
\gamma_{2,5},
\gamma_{2,6},
\gamma_{3,5},
\gamma_{3,6},
\gamma_{4,6}\}.
\end{aligned}
\]

The paper is organized as follows.
In Section~\ref{sec:reduction}, we use Hopf's formula and Pontryagin products to reduce the generating set of $H_2(J_n;\ZZ)$ given in Theorem~\ref{thm:zhang} to $\cB_n\cup\cC_n$.
In Section~\ref{sec:detecting}, we construct a detecting homomorphism $\Theta_n$ from $H_2(J_n;\ZZ)$ to a suitable elementary abelian $2$-group, using a central extension and the embedding of $J_n$ into the semidirect product of the Gauss diagram group and the symmetric group~\cite{BellingeriCheminLebed}.
We also explicitly compute the images of the elements of $\cB_n\cup\cC_n$ under $\Theta_n$.
In Section~\ref{sec:proof}, we prove that $\cB_n\cup\cC_n$ is an $\mathbb F_2$-basis of $H_2(J_n;\ZZ)$ by showing that these images are linearly independent, thereby completing the proof of Theorem~\ref{thm:main}.

\section{Reduction of the generating set}
\label{sec:reduction}

In this section, we use Hopf's formula and Pontryagin products to reduce a generating set of $H_2(J_n;\ZZ)$.
First, we review the definition and basic properties of Pontryagin products.

\subsection{Pontryagin products}
\label{subsec:pontryagin}

Let $G$ be a group, and let $g,h\in G$ be commuting elements.
The Pontryagin product of $g$ and $h$ is defined by
\[
g\wedge h:=
[g|h]-[h|g]
\in H_2(G;\ZZ).
\]
Here $[g|h]-[h|g]$ denotes a $2$-cycle in the normalized bar resolution of $G$.

Now let $G=F/N$, where $F$ is a free group and $N$ is a normal subgroup of $F$.
For commuting elements $g,h\in G$, choose lifts $\widetilde g,\widetilde h\in F$ of $g,h$, respectively.
According to~\cite[Proposition~6.2]{Akita}, the Pontryagin product $g\wedge h\in H_2(G;\ZZ)$ corresponds, under the isomorphism given by Hopf's formula, to
\[
[\widetilde g,\widetilde h]\pmod{[N,F]}.
\]
The right-hand side is independent of the choice of lifts.
Indeed, choose alternative lifts
\[
 \widetilde g'=r\widetilde g,
 \quad
 \widetilde h'=s\widetilde h
 \qquad (r,s\in N).
\]
In the quotient group $F/[N,F]$, the subgroup $N/[N,F]$ is central, so
\[
[r\widetilde g,s\widetilde h]
=[\widetilde g,\widetilde h]
\qquad
\text{in }F/[N,F].
\]
That is,
\[
 [\widetilde g',\widetilde h']
 \equiv[\widetilde g,\widetilde h]
 \pmod{[N,F]}.
\]

In particular, when $G=J_n=F_n/N_n$,
\begin{equation}
\label{eq:pontryagin-hopf}
g\wedge h= \class{[\widetilde g,\widetilde h]}
\in H_2(J_n;\ZZ).
\end{equation}

\begin{proposition}
\label{Prop:Pontryagin}
Let $G$ be a group.
The Pontryagin product on $G$ satisfies the following:
\begin{itemize}
\item[(1)] $g\wedge g=g\wedge 1=1\wedge g=0,$
\item[(2)] $g\wedge h=-h\wedge g$\qquad if $[g,h]=1,$
\item[(3)] 
$\bigl(g_1g_2\bigr)\wedge h=g_1\wedge h+g_2\wedge h$\qquad if $[g_1,h]=[g_2,h]=1,$
\item[(4)]
$g\wedge \bigl(h_1h_2\bigr)=g\wedge h_1+g\wedge h_2$\qquad 
if $[g,h_1]=[g,h_2]=1,$
\item[(5)] 
$(ugu^{-1})\wedge h=g\wedge h$\qquad 
if $[u,h]=[g,h]=1$.
\end{itemize}
\end{proposition}

\begin{proof}
(1), (2), and (4) follow from~\cite[Proposition~6.1]{Akita}.
(3) follows from (2) and (4).

We give a proof of (5).
By assumption, each of $u,g,u^{-1}$ commutes with $h$. Thus (3) gives
\[
(ugu^{-1})\wedge h
 =u\wedge h+g\wedge h+u^{-1}\wedge h.
\]
On the other hand, (1) and (3) give
\[
 0=1\wedge h=(uu^{-1})\wedge h
 =u\wedge h+u^{-1}\wedge h,
\]
and hence (5) follows.
\end{proof}

We next give a proof of Theorem~\ref{thm:zhang} using Hopf's formula and Pontryagin products.

\begin{proof}[Proof of Theorem~\ref{thm:zhang}]
$J_2\cong \ZZ/2$ and
$J_3\cong (\ZZ/2)*(\ZZ/2)$.
Hence, $H_2(J_2;\ZZ)=H_2(J_3;\ZZ)=0$.

Henceforth, let $n\geq 4$ and write $J_n=F_n/N_n$.
For any $r\in N_n$ and $f\in F_n$, we have
\[
frf^{-1}\equiv r\pmod{[N_n,F_n]}.
\]
Moreover, since $[N_n,N_n]\subset[N_n,F_n]$, the group $N_n/[N_n,F_n]$ is abelian and is generated by the images of $A_n\cup B_n\cup C_n$.
Thus each element of $N_n/[N_n,F_n]$ can be represented by an element of $N_n$ of the form
\begin{equation}
\label{eq:zhang-generator}
\prod_{i=2}^n(\sigma_i^2)^{m_i}
\prod_{b(i,j,k)\in B_n}
b(i,j,k)^{n(b(i,j,k))}
\prod_{c(i,j,k)\in C_n}
c(i,j,k)^{n(c(i,j,k))},
\end{equation}
where $m_i$, $n(b(i,j,k))$, and $n(c(i,j,k))$ are integers.

We now consider when the element in~\eqref{eq:zhang-generator} belongs to $[F_n,F_n]$.
The abelianization of the free group $F_n$ is $F_n/[F_n,F_n]\cong \ZZ^{n-1}$, with free basis given by the images of $\sigma_2,\ldots,\sigma_n$.
On the other hand, every element of $B_n$ or $C_n$ is a commutator and hence is trivial in $F_n/[F_n,F_n]$.
Thus the image of~\eqref{eq:zhang-generator} in $F_n/[F_n,F_n]$ is
\[
\prod_{i=2}^n\sigma_i^{2m_i}.
\]
This is trivial in $F_n/[F_n,F_n]$ if and only if
\[
m_2=\cdots=m_n=0.
\]
Therefore, every element of $(N_n\cap[F_n,F_n])/[N_n,F_n]$ can be written in the form
\[
\prod_{b(i,j,k)\in B_n}
b(i,j,k)^{n(b(i,j,k))}
\prod_{c(i,j,k)\in C_n}
c(i,j,k)^{n(c(i,j,k))}.
\]
Consequently, $H_2(J_n;\ZZ)$ is generated by the Hopf classes of the relators in the finite set $B_n\cup C_n$.

We next show that all these Hopf classes are $2$-torsion.
For any $b(i,j,k)=
[\sigma_k\sigma_i\sigma_k,\sigma_j]
\in B_n$,
the elements $\sigma_k\sigma_i\sigma_k$ and $\sigma_j$ commute in $J_n$. Hence~\eqref{eq:pontryagin-hopf} gives
\[
[b(i,j,k)]_{\operatorname H}=
(\sigma_k\sigma_i\sigma_k)\wedge\sigma_j.
\]
Using $\sigma_j^2=1$ and Proposition~\ref{Prop:Pontryagin}, we obtain
\[
2[b(i,j,k)]_{\operatorname H}
=(\sigma_k\sigma_i\sigma_k)\wedge\sigma_j^2
=0.
\]
Thus $[b(i,j,k)]_{\operatorname H}$ is either trivial or has order $2$.
Similarly, $[c(i,j,k)]_{\operatorname H}$ is either trivial or has order $2$.
Therefore, for $n\geq 4$, the group $H_2(J_n;\ZZ)$ is a finite elementary abelian $2$-group.
\end{proof}

\subsection{Reduction of the type \texorpdfstring{$B$}{B} generators}

We first clarify the relationship between the standard generators $\sigma_m=s_{1,m}$ and the general generators $s_{p,q}$.
Applying the interval conjugation relation to $[1,i]\subset[1,k]$ gives
\[
s_{1,k}s_{1,i}s_{1,k}
=s_{k-i+1,k}.
\]
Thus in $J_n$,
\begin{equation}
\label{eq:conjugate-prefix}
\sigma_k\sigma_i\sigma_k=s_{k-i+1,k}.
\end{equation}

\begin{proposition}
\label{prop:B-reduction}
For each
\[
b(i,j,k)\in B_n
\qquad (2\leq i\leq j,\ i+j\leq k\leq n),
\]
we have
\[
\class{b(i,j,k)}
=\class{b(i,j,i+j)}
=\beta_{i,j}.
\]
\end{proposition}

\begin{proof}
Using~\eqref{eq:conjugate-prefix}, define the following elements of $J_n$:
\[
\begin{aligned}
a&:=s_{1,j}=\sigma_j,\\ 
x&:=s_{j+1,i+j}=\sigma_{i+j}\sigma_i\sigma_{i+j},\\
y&:=s_{k-i+1,k}=\sigma_k\sigma_i\sigma_k,\\
u&:=s_{j+1,k}.
\end{aligned}
\]
The index conditions for $b(i,j,k)$ give
\[
 [j+1,i+j]\cap [1,j]=
 [k-i+1,k]\cap [1,j]=
 [j+1,k]\cap [1,j]=\varnothing.
\]
Hence, in $J_n$,
\[
  [x,a]=[y,a]=[u,a]=1.
\]
Moreover, since $[j+1,i+j]\subset[j+1,k]$, the interval conjugation relation gives
\[
\begin{aligned}
  uxu^{-1}
&=s_{j+1,k}s_{j+1,i+j}s_{j+1,k}\\
&=s_{k-i+1,k}=y.
\end{aligned}
\]
By Proposition~\ref{Prop:Pontryagin} (5), we have
\[
y\wedge a
=(uxu^{-1})\wedge a
=x\wedge a.
\]
Here
\begin{align*}
 y\wedge a
&=\class{[\sigma_k\sigma_i\sigma_k,\sigma_j]}
=\class{b(i,j,k)},\\
 x\wedge a
&=\class{[\sigma_{i+j}\sigma_i\sigma_{i+j},\sigma_j]}
=\class{b(i,j,i+j)},
\end{align*}
and hence the assertion follows.
\end{proof}

The following is an immediate consequence of Proposition~\ref{prop:B-reduction}.

\begin{corollary}
\label{cor:B-count}
The subgroup of $H_2(J_n;\ZZ)$ generated by the classes of type $B$ is generated by
\[
 \cB_n=\{\beta_{i,j}\mid 2\leq i\leq j,\ i+j\leq n\}.
\]
\end{corollary}

\begin{remark}
\label{rem:B-count}
For each fixed $i$ with $2\leq i\leq \lfloor n/2\rfloor$, there are $n-2i+1$ choices of $j$, since $i\leq j\leq n-i$.
Therefore,
\[
\#
\{(i,j)\mid 2\leq i\leq j,\ i+j\leq n\}
=\sum_{i=2}^{\lfloor n/2\rfloor}(n-2i+1).
\]
If $n=2m$, then
\[
\begin{aligned}
  \#
\{(i,j)\mid 2\leq i\leq j,\ i+j\leq n\}
  &=\sum_{i=2}^{m}(2m-2i+1)\\
  &=(m-1)^2=\frac{(n-2)^2}{4}.
\end{aligned}
\]
If $n=2m+1$, then
\[
\begin{aligned}
  \#
\{(i,j)\mid 2\leq i\leq j,\ i+j\leq n\}
  &=\sum_{i=2}^{m}(2m+1-2i+1)\\
  &=m(m-1)=\frac{(n-1)(n-3)}{4}.
\end{aligned}
\]
Consequently, we have 
\[
|\cB_n|\leq
\begin{cases}
    \dfrac{(n-2)^2}{4}, & \text{$n$ is even},\\
    \dfrac{(n-1)(n-3)}{4}, & \text{$n$ is odd}.
\end{cases}
\]
\end{remark}

\begin{example}
For $n=6$, we have
\[
\begin{aligned}
  B_6=\{b(2,2,4),b(2,2,5),b(2,2,6),b(2,3,5),b(2,3,6),b(2,4,6),b(3,3,6)\}.  
\end{aligned}
\]
By Proposition~\ref{prop:B-reduction},
\[
\begin{aligned}
 \class{b(2,2,4)}
 &=\class{b(2,2,5)}
=\class{b(2,2,6)}=\beta_{2,2},\\
 \class{b(2,3,5)}
&=\class{b(2,3,6)}=\beta_{2,3},\\
\class{b(2,4,6)}
&=\beta_{2,4},\\
\class{b(3,3,6)}
&=\beta_{3,3}.   
\end{aligned}
\]
Thus the subgroup of $H_2(J_6;\ZZ)$ generated by the classes of type $B$ is generated by
\[
\cB_6=\{
\beta_{2,2},
\beta_{2,3},
\beta_{2,4},
\beta_{3,3}
\}.
\]
\end{example}

\subsection{Reduction of the type \texorpdfstring{$C$}{C} generators}

For a triple $(i,j,k)$ satisfying
\begin{equation}
\label{eq:type-C-condition}
i\geq 1,\  
j\geq 2,\  
i+j<k\leq n,\  
i+j\leq k-i,
\end{equation}
let
\[
A_{i,j,k}:=\sigma_{k-i}\sigma_k\sigma_{i+j}.
\]
Then 
\[
c(i,j,k)=[A_{i,j,k},\sigma_j].
\]

\begin{lemma}
\label{lem:A-factor}
For a triple $(i,j,k)$ satisfying
\eqref{eq:type-C-condition} with $i\geq 2$, we have
\begin{equation}\label{eq:A-factor}
 A_{i,j,k}
 =A_{1,j,k-i+1}A_{i-1,j,k}A_{1,j,i+j}.
\end{equation}
\end{lemma}

\begin{proof}
Using~\eqref{eq:conjugate-prefix}, let
\[
 R:=\sigma_{k-i+1}\sigma_{j+1}\sigma_{k-i+1}
 =s_{k-i-j+1,k-i+1}.
\]
The interval conjugation relation and~\eqref{eq:conjugate-prefix} give
\[
\begin{aligned}
  \sigma_kR\sigma_k
&=s_{1,k} s_{k-i-j+1,k-i+1} s_{1,k}\\
&=s_{i,i+j}\\
&=\sigma_{i+j}\sigma_{j+1}\sigma_{i+j},
\end{aligned}
\]
and hence
\[
R\sigma_k=\sigma_k\sigma_{i+j}\sigma_{j+1}\sigma_{i+j}.
\]
On the other hand, expanding the right-hand side of~\eqref{eq:A-factor} gives
\[
\begin{aligned}
 A_{1,j,k-i+1}A_{i-1,j,k}A_{1,j,i+j}
 &=\sigma_{k-i}R\sigma_k\sigma_{i+j}\sigma_{j+1}\\
 &=\sigma_{k-i}\sigma_k\sigma_{i+j}\sigma_{j+1}\sigma_{i+j}\sigma_{i+j}\sigma_{j+1}\\
 &=\sigma_{k-i}\sigma_k\sigma_{i+j}\\
 &=A_{i,j,k},
\end{aligned}
\]
and hence the assertion follows.
\end{proof}

\begin{proposition}
\label{prop:C-recurrence}
For a triple $(i,j,k)$ satisfying
\eqref{eq:type-C-condition} with $i\geq 2$, we have
\begin{equation}
\label{eq:C-recurrence}
 \class{c(i,j,k)}
 =\class{c(i-1,j,k)}
  +\class{c(1,j,i+j)}
  +\class{c(1,j,k-i+1)}.
\end{equation}
Consequently,
\begin{equation}\label{eq:C-closed}
 \class{c(i,j,k)}
=\sum_{\ell=j+2}^{i+j}\class{c(1,j,\ell)}
 +\sum_{\ell=k-i+1}^{k}\class{c(1,j,\ell)}.
\end{equation}
\end{proposition}

\begin{proof}
First, the triples $(i-1,j,k),\,  (1,j,i+j),\, (1,j,k-i+1)$ satisfy the index conditions for relators of type $C$.
Thus all three factors on the right-hand side of~\eqref{eq:A-factor} commute with $\sigma_j$. 
That is, in $J_n$,
\[
[A_{i-1,j,k},\sigma_j]
=[A_{1,j,i+j},\sigma_j]
=[A_{1,j,k-i+1},\sigma_j]
=1.
\]
Applying Proposition~\ref{Prop:Pontryagin} (3) to~\eqref{eq:A-factor}, we obtain
\[
\begin{aligned}
    \class{c(i,j,k)}
&=A_{i,j,k}\wedge\sigma_j\\
&=(A_{1,j,k-i+1}A_{i-1,j,k}A_{1,j,i+j})\wedge\sigma_j\\
&=A_{1,j,k-i+1}\wedge\sigma_j
+A_{i-1,j,k}\wedge\sigma_j
+A_{1,j,i+j}\wedge\sigma_j\\
&=\class{c(i-1,j,k)}
+\class{c(1,j,i+j)}
+\class{c(1,j,k-i+1)}
\end{aligned}
\]
which proves~\eqref{eq:C-recurrence}.

We next show~\eqref{eq:C-closed} by induction on $i$.
For $i=2$, equation~\eqref{eq:C-recurrence} gives
\[
\begin{aligned}
\class{c(2,j,k)}
&=\class{c(1,j,k)}
+\class{c(1,j,j+2)}
+\class{c(1,j,k-1)}\\
&=\sum_{\ell=j+2}^{j+2}
\class{c(1,j,\ell)}
+\sum_{\ell=k-1}^{k}
\class{c(1,j,\ell)}
\end{aligned}
\]
so~\eqref{eq:C-closed} holds.
Now let $i\geq 3$ and assume that equation~\eqref{eq:C-closed} holds for $i-1$.
Then
\[
\begin{aligned}
\class{c(i,j,k)}
=&\class{c(i-1,j,k)}
+\class{c(1,j,i+j)}
+\class{c(1,j,k-i+1)}\\
=&\sum_{\ell=j+2}^{i+j-1}
\class{c(1,j,\ell)}
+\sum_{\ell=k-i+2}^{k}
\class{c(1,j,\ell)}\\
&\qquad
+\class{c(1,j,i+j)}
+\class{c(1,j,k-i+1)}\\
=&\sum_{\ell=j+2}^{i+j}
\class{c(1,j,\ell)}
 +\sum_{\ell=k-i+1}^{k}
 \class{c(1,j,\ell)}.
\end{aligned}
\]
Thus~\eqref{eq:C-closed} holds by induction.
\end{proof}

The following is an immediate consequence of Proposition~\ref{prop:C-recurrence}.

\begin{corollary}
\label{cor:C-count}
The subgroup of $H_2(J_n;\ZZ)$ generated by the classes of type $C$ is generated by
\[
 \cC_n=
\{\gamma_{j,k}\mid  j\geq 2,\ j+2\leq k\leq n\}.
\]
\end{corollary}

\begin{remark}
\label{rem:C-count}
For each fixed $j$ with $2\leq j\leq n-2$, there are $n-j-1$ choices of $k$, since $j+2\leq k\leq n$.
Therefore,
\[
 \begin{aligned}
\#
\{(j,k)\mid j\geq 2,\ j+2\leq k\leq n\}
&=\sum_{j=2}^{n-2}(n-j-1)\\
&=\frac{(n-3)(n-2)}{2}.
\end{aligned}
\]
Consequently, we have 
\[
 |\cC_n|\leq\frac{(n-3)(n-2)}{2}.
\]
\end{remark}

\begin{example}
For $n=6$, we have
\[
C_6=\{
c(1,2,4),c(1,2,5),c(1,2,6),c(1,3,5),c(1,3,6),c(1,4,6),c(2,2,6)\}.
\]
Applying equation~\eqref{eq:C-closed} to $\class{c(2,2,6)}$ gives
\begin{align*}
\class{c(2,2,6)}
&=
\class{c(1,2,6)}
+\class{c(1,2,4)}
+\class{c(1,2,5)}\\
&=
\gamma_{2,4}
+\gamma_{2,5}
+\gamma_{2,6}.
\end{align*}
Thus the subgroup generated by the classes of type $C$ is generated by
\[
\cC_6=\{
\gamma_{2,4},
\gamma_{2,5},
\gamma_{2,6},
\gamma_{3,5},
\gamma_{3,6},
\gamma_{4,6}
\}.
\]
\end{example}

\section{A homomorphism detecting linear independence}
\label{sec:detecting}

In Section~\ref{sec:reduction}, we obtained a smaller generating set $\cB_n\cup\cC_n$ of $H_2(J_n;\ZZ)$ for $n\geq 4$.
By Remarks~\ref{rem:B-count} and~\ref{rem:C-count}, we have
\[
\begin{aligned}
\dim_{\mathbb F_2}H_2(J_n;\ZZ)
&\leq |\cB_n|+|\cC_n|\\
&\leq
\begin{cases}
\dfrac{(n-2)(3n-8)}{4}, & \text{$n$ is even},\\[2mm]
\dfrac{(n-3)(3n-5)}{4}, & \text{$n$ is odd}.
\end{cases}
\end{aligned}
\]
In this section, we construct a homomorphism $\Theta_n:H_2(J_n;\ZZ)\to Z_n$ to detect the linear independence of the generating set $\cB_n\cup\cC_n$.

\subsection{Codomain of the detecting homomorphism}

Denote the set of all subsets of $\{1,\ldots,n\}$ of cardinality at least $2$ by
\[
\cI_n:=\{I\subset\setn n\mid |I|\geq 2\}.
\]
The symmetric group $S_n$ acts on $\cI_n$ by
\[
\pi(I):=\{\pi(x)\mid x\in I\}
\qquad(\pi\in S_n,\, I\in\cI_n).
\]
Denote the set of all unordered pairs of distinct subsets $I,J\in\cI_n$ by
\[
X_n:=\{\{I,J\}\mid
I,J\in\cI_n,\ I\neq J\}.
\]
The action of the symmetric group $S_n$ on $X_n$ is defined by
\[
\pi\cdot \{I,J\}:=
\{\pi(I),\pi(J)\}\qquad (\pi\in S_n,\, \{I,J\}\in X_n).
\]
For $\{I,J\}\in X_n$, denote its $S_n$-orbit by
\[
\omega(I,J):=
\{\pi\cdot \{I,J\}\mid \pi\in S_n\},
\]
and let the set of orbits be
\[
\cO_n:=
\{\omega(I,J)\mid \{I,J\}\in X_n\}.
\]

For an unordered pair $\{I,J\}\in X_n$, let
\[
a:=\min\{|I|,|J|\},
\qquad
b:=\max\{|I|,|J|\},
\qquad
t:=|I\cap J|.
\]
Then the $S_n$-orbit of $\{I,J\}$ is uniquely determined by the triple $(a,b,t)$.
Such a triple $(a,b,t)$ arises if and only if
\[
2\leq a\leq b\leq n,
\qquad
\max\{0,a+b-n\}\leq t\leq a,
\qquad
(a,b,t)\neq(a,a,a).
\]
Therefore, the orbit set $\cO_n$ is in one-to-one correspondence with the set of all integer triples $(a,b,t)$ satisfying the above conditions.

Let $Z_n$ be the $\mathbb F_2$-vector space with basis $\{z_\omega\mid\omega\in\cO_n\}$.
We also write $z_{a,b,t}$ for the basis element corresponding to the triple $(a,b,t)$.

\subsection{Auxiliary group recording commutator data}

Consider the auxiliary group $\widetilde V_n$ with the following presentation:
\[
\widetilde V_n
=\left\langle
\begin{array}{c}
x_I \quad (I\in\cI_n),\\
z_\omega \quad (\omega\in\cO_n)
\end{array}
\ \middle|\
\begin{array}{l}
x_I^2=1,\quad z_\omega^2=1,\\[1mm]
[z_\omega,x_I]=1,\quad
[z_\omega,z_{\omega'}]=1,\\[1mm]
[x_I,x_J]=z_{\omega(I,J)}
\quad (I\ne J)
\end{array}
\right\rangle.
\]
Here the commutator $[x_I,x_J]$ is recorded by the central element $z_{\omega(I,J)}$.
Moreover, since all commutators in $\widetilde V_n$ are central, for any $a,b,c\in\widetilde V_n$ we have
\begin{equation}
\label{eq:central-commutator}
[ab,c]
=a[b,c]a^{-1}[a,c]
=[a,c][b,c].
\end{equation}

However, the above presentation alone does not make it clear that the central elements $z_\omega$ satisfy no relations other than the defining relations.
The following lemma verifies this and gives a normal form for each element of $\widetilde V_n$.

\begin{lemma}
\label{lem:V-model}
Fix a total ordering $\cI_n=\{I_1<\cdots<I_r\}$.
Then the natural homomorphism
\[
\iota:Z_n
\longrightarrow
\widetilde V_n,
\qquad z_\omega\longmapsto z_\omega
\]
is injective.
Moreover, each element of $\widetilde V_n$ has a unique normal form
\begin{equation}
\label{eq:normal-form}
z\,x_{I_1}^{\varepsilon_1}x_{I_2}^{\varepsilon_2}\cdots x_{I_r}^{\varepsilon_r}\qquad 
(z\in Z_n,\ \varepsilon_a\in \mathbb F_2).
\end{equation}
\end{lemma}

\begin{proof}
First, consider the vector space $W_n$ over $\mathbb F_2$ with standard basis $e_a\,(1\leq a \leq r)$.
Define $c:W_n\times W_n\to Z_n$ by
\[
c(\varepsilon,\delta):=
\sum_{a>b}\varepsilon_a\delta_bz_{\omega(I_a,I_b)}.
\]
The map $c$ is $\mathbb F_2$-linear in each variable and hence is bilinear.
Thus, for any $\varepsilon,\delta,\gamma\in W_n$,
\[
c(\varepsilon,\delta)
+c(\varepsilon+\delta,\gamma)
=c(\delta,\gamma)
+c(\varepsilon,\delta+\gamma).
\]
Therefore, defining a product on the set $Z_n\times W_n$ by
\[
(z,\varepsilon)(z',\delta):=
(z+z'+c(\varepsilon,\delta),\varepsilon+\delta)
\]
gives an associative operation.
The identity element is $(0,0)$, and the inverses are given by
\[
 (z,\varepsilon)^{-1}=(z+c(\varepsilon,\varepsilon),\varepsilon).
\]
Thus, $Z_n\times W_n$ is a group.

Next, consider the elements
\[
\begin{aligned}
  X_a&:=(0,e_a)
  \qquad (a=1,\ldots,r),\\
  Z_\omega&:=(z_\omega,0)
  \qquad (\omega\in\cO_n)
\end{aligned}
\]
of the group $Z_n\times W_n$.
These satisfy all the defining relations of $\widetilde V_n$.
Indeed,
\[
c(e_a,e_a)
=\sum_{1\leq q<p\leq r}(e_a)_p(e_a)_q\,z_{\omega(I_p,I_q)}
=0.
\]
Since $(e_a)_p(e_a)_q\neq 0$ implies $p=q=a$, which contradicts $q<p$, we have $X_a^2=1$.
The definition of $Z_n$ gives $Z_\omega^2=1$.
Moreover, since $c(0,\varepsilon)=c(\varepsilon,0)=0$, each $Z_\omega$ commutes with every $(z,\varepsilon)$.
In particular, $[Z_\omega,X_a]=1,\,
[Z_\omega,Z_{\omega'}]=1$.
For distinct $a,b$, exactly one of $a>b$ and $b>a$ holds, so
\[
c(e_a,e_b)+c(e_b,e_a)
=z_{\omega(I_a,I_b)}.
\]
Hence $X_aX_b=Z_{\omega(I_a,I_b)}X_bX_a$.
Since $X_a^2=X_b^2=1$, we obtain
\[
[X_a,X_b]=Z_{\omega(I_a,I_b)}
\qquad(a\ne b).
\]
Thus the elements $X_a$ and $Z_\omega$ satisfy all the defining relations of $\widetilde V_n$.
The universal property of the presentation therefore gives a homomorphism
\[
\Phi:\widetilde V_n\longrightarrow Z_n\times W_n,
\qquad
x_{I_a}\longmapsto X_a,\quad
z_\omega\longmapsto Z_\omega.
\]
The map $\Phi$ is surjective.
Indeed, for any $(z,\varepsilon)\in Z_n\times W_n$, write
\[
z=\sum_{\omega\in \cO_n}
\delta_\omega z_\omega,
\qquad
\varepsilon=\sum_{a=1}^{r}\varepsilon_a e_a
\qquad
(\delta_\omega,\varepsilon_a\in\mathbb F_2).
\]
By the definition of $c$, for any $1\leq k<r$, we have
\[
\begin{aligned}
c\left(
\sum_{a=1}^{k}\varepsilon_a e_a,
\varepsilon_{k+1}e_{k+1}
\right)
&=\sum_{p>q}
\left(
\sum_{a=1}^{k}\varepsilon_a e_a
\right)_p
(\varepsilon_{k+1}e_{k+1})_q
z_{\omega(I_p,I_q)}\\
&=\sum_{p>k+1}
\left(
\sum_{a=1}^{k}\varepsilon_a e_a
\right)_p
\varepsilon_{k+1}
z_{\omega(I_p,I_{k+1})}\\[2mm]
&=0.
\end{aligned}
\]
By repeatedly applying the definition of the product, we obtain
\[
\begin{aligned}
X_1^{\varepsilon_1}
X_2^{\varepsilon_2}
\cdots X_r^{\varepsilon_r}
&=(0,\varepsilon_1e_1)
  (0,\varepsilon_2e_2)
  \cdots
  (0,\varepsilon_re_r)\\
&=(0,\varepsilon_1e_1+\varepsilon_2e_2) (0,\varepsilon_3e_3)
\cdots
(0,\varepsilon_re_r)\\
&=\cdots\\
&=\left(
0,\sum_{a=1}^{r}\varepsilon_ae_a
\right)
=(0,\varepsilon).
\end{aligned}
\]
Consequently,
\[
\begin{aligned}
(z,\varepsilon)
&=(z,0)(0,\varepsilon)\\
&=
\left(
\prod_{\omega\in \cO_n}
Z_\omega^{\delta_\omega}
\right)
X_1^{\varepsilon_1}
\cdots X_r^{\varepsilon_r}\\
&=
\Phi\left(
\left(
\prod_{\omega\in \cO_n}
z_\omega^{\delta_\omega}
\right)
x_{I_1}^{\varepsilon_1}
\cdots x_{I_r}^{\varepsilon_r}
\right).
\end{aligned}
\]
Thus $\Phi$ is surjective.

We now show that $\iota$ is injective.
Any $z\in Z_n$ can be written uniquely as
\[
z=\sum_{\omega\in \cO_n}
\delta_\omega z_\omega
\qquad
(\delta_\omega\in\mathbb F_2).
\]
Hence
\[
\begin{aligned}
\Phi(\iota(z))
&=\Phi\left(
\prod_{\omega\in \cO_n}
z_\omega^{\delta_\omega}
\right)
=\prod_{\omega\in \cO_n}
\Phi(z_\omega)^{\delta_\omega}
=\prod_{\omega\in \cO_n}
Z_\omega^{\delta_\omega}\\
&=\prod_{\omega\in \cO_n}
(\delta_\omega z_\omega,0)
=\left(
\sum_{\omega\in \cO_n}
\delta_\omega z_\omega,0
\right)
=(z,0).
\end{aligned}
\]
For any $z\in\Ker\iota$, we have $\iota(z)=1$.
Applying $\Phi$ to both sides gives $(z,0)=(0,0)$, so $z=0$.
Thus $\iota$ is injective, and $Z_n$ can be identified with its image $\iota(Z_n)$.

Finally, consider normal forms for elements of $\widetilde V_n$. 
By repeatedly applying the defining relations of $\widetilde V_n$, any $g\in\widetilde V_n$ can be written as
\[
g=
\left(
\prod_{\omega\in \cO_n}
z_\omega^{\delta_\omega}
\right)
x_{I_1}^{\varepsilon_1}
x_{I_2}^{\varepsilon_2}
\cdots x_{I_r}^{\varepsilon_r}
\qquad
(\delta_\omega,\varepsilon_a\in\mathbb F_2).
\]
Let
\[
z:=
\sum_{\omega\in \cO_n}
\delta_\omega z_\omega\in Z_n.
\]
By the injectivity of $\iota$,
\[
\iota(z)= \prod_{\omega\in \cO_n}
z_\omega^{\delta_\omega}\in \widetilde V_n
\]
can be identified with $z$, giving the normal form~\eqref{eq:normal-form}.
Now suppose that
\[
z\,x_{I_1}^{\varepsilon_1}\cdots x_{I_r}^{\varepsilon_r}
=z'\,x_{I_1}^{\varepsilon_1'}\cdots x_{I_r}^{\varepsilon_r'}.
\]
Applying $\Phi$ to both sides gives
\[
\left(
z,\sum_{a=1}^{r}\varepsilon_a e_a
\right)
=\left(
z',\sum_{a=1}^{r}\varepsilon_a'e_a
\right),
\]
which implies
\[
z=z',
\qquad
\varepsilon_a=\varepsilon_a'\quad(1\leq a\leq r).
\]
Therefore, the normal form~\eqref{eq:normal-form} is unique.
\end{proof}

\subsection{Central extension of a semidirect product}
\label{subsec:central-extension}

Let $V_n$ be the $\mathbb F_2$-vector space with basis $\{e_I\mid I\in\cI_n\}$.
In this subsection, we construct the following central extension:
\begin{equation}
\label{eq:central-extension}
1\longrightarrow Z_n
\longrightarrow
\widetilde V_n\rtimes S_n
\longrightarrow
V_n\rtimes S_n
\longrightarrow 1.
\end{equation}

\begin{lemma}
\label{lem:V-central-extension}
Define a map $q:\widetilde V_n\to V_n$ by
\[
 x_I\longmapsto e_I,
 \qquad z_\omega\longmapsto 0.
\]
Then we obtain the following central extension:
\begin{equation}
\label{eq:V-central-extension}
1\longrightarrow Z_n
\xrightarrow{ \ \iota \ }
\widetilde V_n
\xrightarrow{ \ q \ }
V_n
\longrightarrow 1.
\end{equation}
\end{lemma}

\begin{proof}
First, since the assignment $x_I\mapsto e_I,\,z_\omega\mapsto 0$ satisfies all the defining relations of $\widetilde V_n$, the universal property of the presentation shows that $q$ is a group homomorphism.
Moreover, each basis element $e_I$ of $V_n$ is the image of $x_I$, so $q$ is surjective.
By Lemma~\ref{lem:V-model}, the natural homomorphism $\iota:Z_n\to\widetilde V_n$ is injective, and we identify $Z_n$ with its image $\iota(Z_n)$.

Next, we have $\Ker q=Z_n$.
Indeed, by the normal form~\eqref{eq:normal-form}, any $g\in\Ker q$ can be written uniquely as
\[
g=
z\,x_{I_1}^{\varepsilon_1}
x_{I_2}^{\varepsilon_2}
\cdots
x_{I_r}^{\varepsilon_r}
\qquad
(z\in Z_n,\,
\varepsilon_a\in\mathbb F_2).
\]
Hence
\[
0=q(g)
=q\left(
z\,x_{I_1}^{\varepsilon_1}
\cdots
x_{I_r}^{\varepsilon_r}
\right)
=\sum_{a=1}^{r}\varepsilon_a e_{I_a}.
\]
Since $e_{I_a}\,(1\leq a \leq r)$ form a basis of $V_n$, we have $\varepsilon_a=0\,(1\leq a \leq r)$ and hence $g=z\in Z_n$.
Thus $\Ker q\subset Z_n$. 
On the other hand, since $Z_n$ is generated by $z_\omega\, (\omega\in \cO_n)$, any $z\in Z_n$ can be written as
\[
z=
\prod_{\omega\in \cO_n}
z_\omega^{\delta_\omega}
\qquad
(\delta_\omega\in\mathbb F_2).
\]
Since $q$ is a homomorphism and $q(z_\omega)=0$ for each $\omega\in \cO_n$, we obtain
\[
q(z)
=q\left(
\prod_{\omega\in \cO_n}
z_\omega^{\delta_\omega}
\right)
=\sum_{\omega\in \cO_n}
\delta_\omega q(z_\omega)
=0.
\]
Hence $z\in\Ker q$, which also gives $Z_n\subset \Ker q$.
Therefore, \eqref{eq:V-central-extension} is a short exact sequence.

Finally, the basis elements $z_\omega$ of $Z_n$ lie in the center of $\widetilde V_n$, so $Z_n\subset Z(\widetilde V_n)$ and~\eqref{eq:V-central-extension} is a central extension.
\end{proof}

We next construct the semidirect product $\widetilde V_n\rtimes S_n$.
First, for each permutation $\pi\in S_n$, the elements $x_{\pi(I)},z_\omega$ satisfy all the defining relations of $\widetilde V_n$.
In particular, for $I\neq J$,
\[
[x_{\pi(I)},x_{\pi(J)}]
=z_{\omega(\pi(I),\pi(J))}
=z_{\omega(I,J)}.
\]
The last equality holds because $\omega(I,J)$ is the $S_n$-orbit of the unordered pair $\{I,J\}$.
Thus, the universal property of the presentation gives a homomorphism
\[
\alpha_\pi:
\widetilde V_n\longrightarrow\widetilde V_n,\qquad
x_I\longmapsto x_{\pi(I)},
\quad
z_\omega\longmapsto z_\omega.
\]
Similarly, $\pi^{-1}\in S_n$ gives a homomorphism
\[
\alpha_{\pi^{-1}}:
\widetilde V_n\longrightarrow\widetilde V_n,\qquad 
x_I \longmapsto x_{\pi^{-1}(I)},\quad
z_\omega \longmapsto z_\omega.
\]
Since their compositions fix every generator $x_I,z_\omega$, we have $\alpha_\pi^{-1}=\alpha_{\pi^{-1}}$.
Thus $\alpha_\pi$ is an automorphism of $\widetilde V_n$.
Moreover, for any $\pi,\,\rho\in S_n$, we have $\alpha_\pi\circ\alpha_\rho=\alpha_{\pi\rho}$, giving a group homomorphism
\[
\varphi:S_n\longrightarrow
\operatorname{Aut}(\widetilde V_n),
\qquad
\pi\longmapsto\alpha_\pi.
\]
We define the semidirect product $\widetilde V_n\rtimes_\varphi S_n$ with respect to the action $\varphi$.
Its multiplication is given by
\[
(g,\pi)(h,\rho)
=(g\alpha_\pi(h),\pi\rho).
\]

Similarly, for any $\pi\in S_n$, the map
\[
\beta_\pi:V_n
\longrightarrow V_n,\qquad e_I\longmapsto e_{\pi(I)}
\]
is an automorphism of $V_n$, and these maps define a group homomorphism
\[
\psi:S_n \longrightarrow\operatorname{Aut}(V_n),\qquad
\pi \longmapsto \beta_\pi.
\]
Therefore, we obtain the semidirect product $V_n\rtimes_\psi S_n$.

The surjection $q:\widetilde V_n\to V_n$ is $S_n$-equivariant with respect to the actions $\alpha_\pi$ and $\beta_\pi$.
That is, for any $\pi\in S_n$,
\begin{equation}
\label{eq:equivariant-alpha-beta}
q\circ\alpha_\pi
=\beta_\pi\circ q.
\end{equation}

\begin{proposition}
\label{prop:central-extension}
Define the map between semidirect products
\[
q\rtimes\operatorname{id}_{S_n}:
\widetilde V_n\rtimes_\varphi S_n
\longrightarrow
V_n\rtimes_\psi S_n,\qquad
(g,\pi)
\longmapsto(q(g),\pi)
\]
and the homomorphism
\[
\iota':
Z_n\longrightarrow
\widetilde V_n\rtimes_\varphi S_n,
\qquad
z\longmapsto(z,1_{S_n}).
\]
Then we obtain the central extension~\eqref{eq:central-extension}.
\end{proposition}

\begin{proof}
First, for any
$g,h\in\widetilde V_n$ and
$\pi,\rho\in S_n$,
\[
\begin{aligned}
(q\rtimes\operatorname{id}_{S_n})
((g,\pi)(h,\rho))
&=(q\rtimes\operatorname{id}_{S_n})
(g\alpha_\pi(h),\pi\rho)\\
&=(q(g)+q(\alpha_\pi(h)),\pi\rho),\\
(q\rtimes\operatorname{id}_{S_n})(g,\pi)(q\rtimes\operatorname{id}_{S_n})(h,\rho)
&=(q(g),\pi)(q(h),\rho)\\
&=(q(g)+\beta_\pi(q(h)),\pi\rho).
\end{aligned}
\]
Then $q\rtimes\operatorname{id}_{S_n}$ is a homomorphism by~\eqref{eq:equivariant-alpha-beta}.
By the surjectivity of $q$, for any $(v,\pi)\in V_n\rtimes_\psi S_n$, there exists $g\in\widetilde V_n$ such that $q(g)=v$.
Then
\[
(q\rtimes\operatorname{id}_{S_n})(g,\pi)
=(v,\pi).
\]
Hence $q\rtimes\operatorname{id}_{S_n}$ is surjective.
Moreover, the injectivity of $\iota$ implies that $\iota'$ is also injective, and we identify $Z_n$ with its image $\iota'(Z_n)$.

We next compute $\Ker(q\rtimes\operatorname{id}_{S_n})$.
For any $(g,\pi)\in\widetilde V_n\rtimes_\varphi S_n$,
\[
\begin{aligned}
(g,\pi)\in\Ker(q\rtimes\operatorname{id}_{S_n})
&\Longleftrightarrow
(q(g),\pi)=(0,1_{S_n})\\
&\Longleftrightarrow
g\in\Ker q=Z_n,\ \pi=1_{S_n}.
\end{aligned}
\]
Hence
\[
\Ker(q\rtimes\operatorname{id}_{S_n})
=\left\{
(z,1_{S_n})\mid z\in Z_n
\right\}
\cong Z_n.
\]
Thus \eqref{eq:central-extension} is a short exact sequence.

Finally, for any $z\in Z_n$ and $(g,\pi)\in\widetilde V_n\rtimes_\varphi S_n$,
we have
\[
\begin{aligned}
(z,1_{S_n})(g,\pi)
&=(zg,\pi)
=(gz,\pi)\\
&=(g\alpha_\pi(z),\pi)
=(g,\pi)(z,1_{S_n}).
\end{aligned}
\]
Therefore, $Z_n$ lies in the center of $\widetilde V_n\rtimes_\varphi S_n$ and~\eqref{eq:central-extension} is a central extension.
\end{proof}

\subsection{Construction of the detecting homomorphism}

Consider the Gauss diagram group $D_n$ defined by the following presentation (see~\cite[Section~1]{BellingeriCheminLebed}):
\[
D_n=
\left\langle
\tau_I\ (I\in\cI_n)
\ \middle|\
\begin{array}{l}
\tau_I^2=1,\\[1mm]
[\tau_I,\tau_J]=1
\quad
(I\subset J \text{ or } I\cap J=\varnothing)
\end{array}
\right\rangle.
\]
For any $\pi\in S_n$, the elements $\tau_{\pi(I)}$ satisfy all the defining relations of $D_n$, so the universal property of the presentation gives a homomorphism
\[
\delta_\pi:D_n\longrightarrow D_n,\qquad
\tau_I\longmapsto\tau_{\pi(I)}.
\]
Moreover, the homomorphism
\[
\delta_{\pi^{-1}}:D_n\longrightarrow D_n,
\qquad
\tau_I\longmapsto\tau_{\pi^{-1}(I)}
\]
is the inverse of $\delta_\pi$. 
Thus $\delta_\pi$ is an automorphism of $D_n$.
For any $\pi,\rho\in S_n$, we have 
$\delta_{\pi\rho}
=\delta_\pi\circ\delta_\rho$, which gives a group homomorphism
\[
\lambda:S_n
\longrightarrow\operatorname{Aut}(D_n),\qquad
\pi
\longmapsto \delta_\pi.
\]
Thus we define the semidirect product $D_n\rtimes_\lambda S_n$ with respect to the action $\lambda$.

There is also a natural surjective homomorphism
\[
s:J_n\longrightarrow S_n,
\qquad
s_{p,q}\longmapsto r_{p,q},
\]
where $r_{p,q}$ denotes the reversal of $[p,q]$.

As explained in~\cite[Section~1]{BellingeriCheminLebed}, there is an injective homomorphism
\[
\rho:
J_n
\longrightarrow
D_n\rtimes_\lambda S_n,\qquad
\rho(s_{p,q})=(\tau_{[p,q]},r_{p,q}).
\]
In particular, since $\sigma_m=s_{1,m}$, we have
\begin{equation}
\label{eq:rho-sigma}
\rho(\sigma_m)
=(\tau_{K_m},r_m),
\end{equation}
where
\[
K_m:=\{1,\ldots,m\},
\qquad
r_m:=r_{1,m}.
\]

The abelianization of $D_n$ is $D_n^{\operatorname{ab}}\cong V_n$, with abelianization map
\[
\operatorname{ab}:D_n\longrightarrow V_n,
\qquad
\tau_I\longmapsto e_I.
\]
The map $\operatorname{ab}$ is $S_n$-equivariant with respect to the actions $\delta_\pi$ and $\beta_\pi$.
That is, for any $\pi\in S_n$, we have $\operatorname{ab}\circ\delta_\pi
=\beta_\pi\circ\operatorname{ab}$.
Thus we obtain a homomorphism between semidirect products:
\[
\begin{aligned}
\operatorname{ab}\rtimes\operatorname{id}_{S_n}:
D_n\rtimes_\lambda S_n
&\longrightarrow
V_n\rtimes_\psi S_n,\\
(g,\pi)
&\longmapsto
(\operatorname{ab}(g),\pi).
\end{aligned}
\]
Define
\[
\overline\rho
:=
(\operatorname{ab}\rtimes\operatorname{id}_{S_n})
\circ\rho.
\]
Then the following diagram commutes:
\[
\begin{tikzcd}[column sep=large,row sep=large]
J_n
\arrow[r,"\rho"]
\arrow[dr,"\overline\rho"']
&
D_n\rtimes_\lambda S_n
\arrow[d,"\operatorname{ab}\rtimes\operatorname{id}_{S_n}"]
\\
&
V_n\rtimes_\psi S_n.
\end{tikzcd}
\]
Moreover, by~\eqref{eq:rho-sigma}, we have
\[
\overline\rho(\sigma_m)
=(e_{K_m},r_m).
\]

Furthermore, define a homomorphism on the free group $F_n=F(\sigma_2,\ldots,\sigma_n)$ by
\begin{equation}\label{eq:free-lift}
\widetilde\rho:F_n
\longrightarrow
\widetilde V_n\rtimes_\varphi S_n,\qquad
\sigma_m
\longmapsto
(x_{K_m},r_m).
\end{equation}
Then the following diagram commutes:
\[
\begin{tikzcd}
F_n
\arrow[r,"\widetilde\rho"]
\arrow[d,"p"']
&
\widetilde V_n\rtimes_{\varphi}S_n
\arrow[d,"q\rtimes\operatorname{id}_{S_n}"]
\\
J_n
\arrow[r,"\overline\rho"']
&
V_n\rtimes_{\psi}S_n.
\end{tikzcd}
\]
Indeed, for each free generator $\sigma_m\in F_n$,
\[
\begin{aligned}
(q\rtimes\operatorname{id}_{S_n})
(\widetilde\rho(\sigma_m))
&=
(q\rtimes\operatorname{id}_{S_n})
(x_{K_m},r_m)\\
&=
(e_{K_m},r_m)\\
&=
\overline\rho(p(\sigma_m)),
\end{aligned}
\]
and hence
$(q\rtimes\operatorname{id}_{S_n})
\circ\widetilde\rho
=\overline\rho\circ p$.

The above constructions can be summarized by the following commutative diagram:
\[
\begin{tikzcd}[column sep=large,row sep=large]
F_n
\arrow[r,"p"]
\arrow[d,"\widetilde\rho"']
&
J_n
\arrow[r,"\rho"]
\arrow[d,"\overline\rho"']
&
D_n\rtimes_\lambda S_n
\arrow[d,"\operatorname{ab}\rtimes\operatorname{id}_{S_n}"]
\\
\widetilde V_n\rtimes_\varphi S_n
\arrow[r,"q\rtimes\operatorname{id}_{S_n}"']
&
V_n\rtimes_\psi S_n
\arrow[r,equal]
&
V_n\rtimes_\psi S_n.
\end{tikzcd}
\]

We now construct a homomorphism to detect the linear independence of $\cB_n\cup\cC_n$ in $H_2(J_n;\ZZ)$.

\begin{proposition}
\label{prop:Theta-construction}
$\widetilde\rho$ induces a homomorphism
\begin{equation}
\label{eq:Theta}
\Theta_n:
\frac{N_n\cap[F_n,F_n]}{[N_n,F_n]}
\longrightarrow
Z_n,
\qquad
\class{r}\longmapsto\widetilde\rho(r).
\end{equation}
By Hopf's formula, this can be regarded as
\[
\Theta_n:H_2(J_n;\ZZ)\longrightarrow Z_n.
\]
\end{proposition}

\begin{proof}
Since $p(r)=1$ for any $r\in N_n$, we have
\[
(q\rtimes\operatorname{id}_{S_n})(\widetilde\rho(r))
=(\overline\rho\circ p)(r)
=\overline\rho(1)
=(0,1_{S_n}).
\]
Thus
\[
\widetilde\rho(r)
\in\Ker(q\rtimes\operatorname{id}_{S_n})
\cong Z_n.
\]
Therefore, the restriction of $\widetilde\rho$ to $N_n$ defines a homomorphism
\[
\widetilde\rho|_{N_n}:
N_n\longrightarrow Z_n,\qquad
r \longmapsto \widetilde\rho(r).
\]

For any $r\in N_n$ and $f\in F_n$, since $\widetilde\rho(r)\in Z_n$ and $Z_n\subset Z(\widetilde V_n\rtimes_\varphi S_n)$, we have
\[
\widetilde\rho([r,f])
=[\widetilde\rho(r),\widetilde\rho(f)]
=1.
\]
Hence
\[
[N_n,F_n]\subset
\Ker(\widetilde\rho|_{N_n}).
\]
Thus $\widetilde\rho|_{N_n}$ induces a unique homomorphism
\[
\widetilde\rho':
\frac{N_n}{[N_n,F_n]}
\longrightarrow Z_n,
\qquad
r\pmod{[N_n,F_n]}
\longmapsto
\widetilde\rho(r).
\]

The inclusions
\[
[N_n,F_n]\subset N_n\cap[F_n,F_n]
\subset N_n
\]
induce a natural inclusion
\[
\frac{N_n\cap[F_n,F_n]}{[N_n,F_n]}
\hookrightarrow
\frac{N_n}{[N_n,F_n]}.
\]
By restricting $\widetilde\rho'$ to this subgroup, we obtain the homomorphism~\eqref{eq:Theta}.

The construction of $\Theta_n$ can be summarized by the following commutative diagram:
\[
\begin{tikzcd}[column sep=large,row sep=large]
N_n
\arrow[r,two heads,""]
\arrow[d,"\widetilde\rho|_{N_n}"']
&
\dfrac{N_n}{[N_n,F_n]}
\arrow[d,"\widetilde\rho'"]
&
\dfrac{N_n\cap[F_n,F_n]}{[N_n,F_n]}
\arrow[l,hook]
\arrow[d,"\Theta_n"]
\\
Z_n
\arrow[r,equal]
&
Z_n
&
Z_n
\arrow[l,equal].
\end{tikzcd}
\]
\end{proof}

\subsection{Images of the generators under the detecting homomorphism}

\begin{lemma}
\label{lem:semidirect-commutator}
Let a group $S$ act on a group $V$ by automorphisms, and consider the semidirect product $V\rtimes S$.
Let $g=(v,p)$ and $h=(w,q)$, and suppose that
\[
pq=qp,\qquad p(w)=w,\qquad 
q(v)=v.
\]
Then
\[
[g,h]=([v,w],1_S)\in V\rtimes S.
\]
\end{lemma}

\begin{proof}
By the multiplication rule in the semidirect product and the assumptions, we have
\[
gh=(vw,pq),\qquad hg=(wv,pq).
\]
Hence
\[
([v,w],1_S)(hg)
=([v,w]wv,pq)
=(vw,pq)
=gh,
\]
and the assertion follows.
\end{proof}

\begin{proposition}
\label{prop:detector-images}
The homomorphism $\Theta_n$ satisfies the following:
\begin{alignat}{2}
\Theta_n(\beta_{i,j})
&= z_{i,j,0}
&\qquad&
(2\leq i\leq j,\ i+j\leq n),
\label{eq:beta-image}\\
\Theta_n(\gamma_{j,k})
&= z_{j,k-1,j}+z_{j,k,j}+z_{j,j+1,j}
&&
(j\geq 2,\ j+2\leq k\leq n).
\label{eq:gamma-image}
\end{alignat}
\end{proposition}

\begin{proof}
We first compute $\Theta_n(\beta_{i,j})$. 
Let
\[
\begin{aligned}
I&:=\{1,\ldots,i\},&\qquad
X&:=\{1,\ldots,j\},\\
K&:=\{1,\ldots,i+j\},&\qquad
Y&:=\{j+1,\ldots,i+j\},
\end{aligned}
\]
so that $Y=r_{i+j}(I)$. By~\eqref{eq:free-lift}, we obtain
\[
\begin{aligned}
\widetilde\rho(\sigma_{i+j}\sigma_i\sigma_{i+j})
&=(x_K,r_{i+j})(x_I,r_i)(x_K,r_{i+j})\\
&=(x_Kx_{r_{i+j}(I)}x_{r_{i+j}r_i(K)},
r_{i+j}r_ir_{i+j})\\
&=(x_Kx_Yx_K,r_Y),
\end{aligned}
\]
where $r_Y$ denotes the permutation reversing the interval $Y$.
The last equality follows from $r_{i+j}(I)=Y$, $r_i(K)=r_{i+j}(K)=K$, and $r_{i+j}r_i r_{i+j}=r_Y$.
We also have
\[
\widetilde\rho(\sigma_j)=(x_X,r_X).
\]

Since $X\cap Y=\varnothing$ and $K=X\cup Y$,
\[
r_Xr_Y=r_Yr_X,\qquad
r_Y(x_X)=x_X,\qquad
r_X(x_Kx_Yx_K)=x_Kx_Yx_K.
\]
Therefore, identifying $Z_n$ with its image $\iota'(Z_n)\subset\widetilde V_n\rtimes_\varphi S_n$, we obtain from Lemma~\ref{lem:semidirect-commutator} that
\[
\begin{aligned}
\Theta_n(\beta_{i,j})
&=\widetilde\rho(
[\sigma_{i+j}\sigma_i\sigma_{i+j},\sigma_j])\\
&=[(x_Kx_Yx_K,r_Y),(x_X,r_X)]\\
&=([x_Kx_Yx_K,x_X],1_{S_n})\\
&=[x_Kx_Yx_K,x_X].
\end{aligned}
\]
By~\eqref{eq:central-commutator}, we further obtain
\[
\begin{aligned}
\Theta_n(\beta_{i,j})
&=[x_K,x_X]+[x_Y,x_X]+[x_K,x_X]\\
&=[x_Y,x_X]
=z_{\omega(X,Y)}
=z_{i,j,0},
\end{aligned}
\]
which proves~\eqref{eq:beta-image}.

We next compute $\Theta_n(\gamma_{j,k})$. 
Let
\[
\begin{aligned}
X&:=\{1,\ldots,j\},&\qquad
I_1&:=\{1,\ldots,k-1\},\\
I_2&:=\{1,\ldots,k\},&\qquad
I_3&:=X\cup\{k\}.
\end{aligned}
\]
By~\eqref{eq:free-lift}, we have
\[
\begin{aligned}
  \widetilde\rho(\sigma_{k-1}\sigma_k\sigma_{j+1})
&=(x_{I_1},r_{k-1})(x_{I_2},r_k)
(x_{K_{j+1}},r_{j+1})\\
&=(x_{I_1}x_{r_{k-1}(I_2)}
x_{r_{k-1}r_k(K_{j+1})},
r_{k-1}r_kr_{j+1})\\
&=(x_{I_1}x_{I_2}x_{I_3},
r_{k-1}r_kr_{j+1}).  
\end{aligned}
\]
The last equality follows from $r_{k-1}(I_2)=I_2$ and $r_{k-1}r_k(K_{j+1})=I_3$.

Let $p:=r_{k-1}r_kr_{j+1}$.
For any $t\in X$, we have
\[
p(t)=j+1-t=r_X(t),
\]
so $p(X)=X$ and $p|_X=r_X|_X$.
Since $p$ preserves $X^c$ and $r_X$ acts as the identity on $X^c$, we obtain
\[
pr_X=r_Xp.
\]
Moreover, $r_X$ preserves each of $I_1,I_2,I_3$ as a set, and hence
\[
r_X(x_{I_1}x_{I_2}x_{I_3})
=x_{I_1}x_{I_2}x_{I_3}.
\]
Since $p(X)=X$, we also have
\[
p(x_X)=x_X.
\]
Therefore, identifying $Z_n$ with its image $\iota'(Z_n)$, we obtain by Lemma~\ref{lem:semidirect-commutator} that
\[
\begin{aligned}
\Theta_n(\gamma_{j,k})
&=\widetilde\rho(
[\sigma_{k-1}\sigma_k\sigma_{j+1},\sigma_j])\\
&=[
(x_{I_1}x_{I_2}x_{I_3},r_{k-1}r_kr_{j+1}),
(x_X,r_X)]\\
&=([x_{I_1}x_{I_2}x_{I_3},x_X],1_{S_n})\\
&=[x_{I_1}x_{I_2}x_{I_3},x_X].
\end{aligned}
\]
Applying~\eqref{eq:central-commutator}, we have
\[
\begin{aligned}
\Theta_n(\gamma_{j,k})
&=[x_{I_1},x_X]+[x_{I_2},x_X]+[x_{I_3},x_X]\\
&=z_{j,k-1,j}+z_{j,k,j}+z_{j,j+1,j},
\end{aligned}
\]
which proves~\eqref{eq:gamma-image}.
\end{proof}

\section{Proof of the main theorem}
\label{sec:proof}

\begin{proof}[Proof of Theorem~\ref{thm:main}]

Let $n\geq 4$.
It suffices to show that the generating set
$\cB_n\cup\cC_n$ is linearly independent
over $\mathbb F_2$ in $H_2(J_n;\ZZ)$.

Let $b_{i,j},\,c_{j,k}\in\mathbb F_2$ be coefficients satisfying
\[
\sum_{\substack{2\leq i\leq j\\i+j\leq n}}
b_{i,j}\beta_{i,j}
+
\sum_{\substack{j\geq2\\j+2\leq k\leq n}}
c_{j,k}\gamma_{j,k}=0.
\]
Applying $\Theta_n$ to both sides and using Proposition~\ref{prop:detector-images}, we obtain
\begin{equation}
\label{eq:independence-relation}
0=
\sum_{\substack{2\leq i\leq j\\i+j\leq n}}
b_{i,j}z_{i,j,0}
+\sum_{\substack{j\geq 2\\j+2\leq k\leq n}}
c_{j,k}(z_{j,k-1,j}+z_{j,k,j}+z_{j,j+1,j}).
\end{equation}

Since the basis elements $z_{i,j,0}$ are distinct for distinct pairs $(i,j)$ with $i\leq j$ and appear only in the first sum, comparing their coefficients in \eqref{eq:independence-relation} yields
\[
b_{i,j}=0\qquad
(2\leq i\leq j,\ i+j\leq n).
\]

For each $j$ with $2\leq j\leq n-2$, let
\[
\eta_{j,t}:=z_{j,t,j}
\qquad (j+1\leq t\leq n).
\]
The sets $\{\eta_{j,t}\mid j+1\leq t\leq n\}$ for $2\leq j\leq n-2$ are pairwise disjoint subsets of the basis of $Z_n$.
Hence, for each fixed $j$, \eqref{eq:independence-relation} gives
\[
\begin{aligned}
0
&=\sum_{k=j+2}^{n}
c_{j,k}
(\eta_{j,k-1}+\eta_{j,k}+\eta_{j,j+1})\\
&=\left(\sum_{k=j+3}^{n}c_{j,k}\right)
\eta_{j,j+1}
+\sum_{t=j+2}^{n-1}
(c_{j,t}+c_{j,t+1})\eta_{j,t}
+c_{j,n}\eta_{j,n}.
\end{aligned}
\]
Comparing the coefficient of $\eta_{j,n}$, we obtain
\[
c_{j,n}=0.
\]
Furthermore, comparing the coefficients of $\eta_{j,t}$ for $t=n-1,n-2,\ldots,j+2$ gives
\[
c_{j,t}+c_{j,t+1}=0.
\]
Therefore, by downward induction, we have
\[
c_{j,k}=0
\qquad (j+2\leq k\leq n).
\]

Consequently, all coefficients $b_{i,j}$ and $c_{j,k}$ vanish, and hence $\cB_n\cup\cC_n$ is linearly independent in $H_2(J_n;\ZZ)$.
This proves Theorem~\ref{thm:main}.
\end{proof}

\begin{ack}
The author would like to thank Toshiyuki Akita for valuable comments and suggestions. 
The author also gratefully acknowledges financial support from the Department of Mathematics at Hokkaido University.
\end{ack}

\end{document}